\documentclass[final,onefignum,onetabnum]{siamart190516}

\usepackage{subcaption}
\usepackage{lipsum}
\usepackage{amsfonts}
\usepackage{amssymb}
\usepackage{graphicx}
\usepackage{epstopdf}
\usepackage{algorithmic}
\ifpdf
  \DeclareGraphicsExtensions{.eps,.pdf,.png,.jpg}
\else
  \DeclareGraphicsExtensions{.eps}
\fi

\newsiamremark{remark}{Remark}
\newsiamremark{hypothesis}{Hypothesis}
\crefname{hypothesis}{Hypothesis}{Hypotheses}
\newsiamthm{claim}{Claim}

\headers{Theory for Accelerated Optimization}{I. M. Ross}

\title{An Optimal Control Theory for Accelerated Optimization
}

\author{Isaac M. Ross\thanks{Distinguished Professor, Naval Postgraduate School, Monterey, CA
  (\email{imross@nps.edu}}).}

\usepackage{amsopn}

\newcommand{\Real}{\mathbb R}
\newcommand{\set}[1]{\left\{#1\right\}}
\newcommand{\real}[1]{{\mathbb R}^{#1}}

\newcommand{\bb}{{\boldsymbol b}}
\newcommand{\bc}{{\boldsymbol c}}

\newcommand{\be}{{\boldsymbol e}}
\newcommand{\bff}{{\boldsymbol f}}

\newcommand{\bs}{{\boldsymbol s}}
\newcommand{\bu}{{\boldsymbol u}}

\newcommand{\bv}{{\boldsymbol v}}

\newcommand{\bx}{\boldsymbol x}

\newcommand{\bz}{{\boldsymbol z}}

\newcommand{\buf}{{\bu(\cdot)}}  

\newcommand{\bM}{{\boldsymbol M}}

\newcommand{\bQ}{{\boldsymbol Q}}

\newcommand{\bW}{{\boldsymbol W}}

\newcommand{\U}{\mathbb{U}}

\newcommand{\bzero}{{\bf 0}}

\newcommand{\bnu}{{\mbox{\boldmath $\nu$}}}
\newcommand{\bmu}{{\mbox{\boldmath $\mu$}}}
\newcommand{\bpsi}{{\mbox{\boldmath $\psi$}}}
\newcommand{\blam}{{\mbox{\boldmath $\lambda$}}}

\newcommand{\bchi}{\mbox{\boldmath$\chi$}}
\newcommand{\bomega}{\mbox{\boldmath$\omega$}}

\newcommand{\btheta}{\mbox{\boldmath$\theta$}}
\newcommand{\bzeta}{\mbox{\boldmath$\zeta$}}

\ifpdf
\hypersetup{
  pdftitle={An Optimal Control Theory for Accelerated Optimization},
  pdfauthor={I. M. Ross}
}
\fi

\usepackage[utf8]{inputenc}

\begin{document}

\maketitle

\begin{abstract}
The first-order optimality conditions for a generic nonlinear optimization problem are generated as part of the terminal transversality conditions of an optimal control problem. It is shown that the Lagrangian of the optimization problem is connected to the Hamiltonian of the optimal control problem via a zero-Hamiltonian, infinite-order, singular arc. The necessary conditions for the singular optimal control problem are used to produce an auxiliary controllable dynamical system whose trajectories generate algorithm primitives for the optimization problem.  A three-step iterative map for a generic algorithm is designed by a semi-discretization step. Neither the feedback control law nor the differential equation governing the algorithm need be derived explicitly. A search direction is produced by a proximal-aiming-type method that dissipates a control Lyapunov function.  New step size procedures based on minimizing control Lyapunov functions along a search vector complete the design of the accelerated algorithms.

\end{abstract}



\section{Introduction}
Consider a generic, nonlinear optimization problem,
\begin{equation}\label{eq:NLP}
(N) \left\{\displaystyle\mathop\text{Minimize }_{\bx_f \in C \subseteq \real{N_x}}   E(\bx_f) \right.
\end{equation}
where, $E: \real{N_x} \ni \bx_f \mapsto \Real$ is an objective function,  $C$ is a constraint set in $\real{N_x}$ and $ N_x \in \mathbb{N}^+$.  The first-order optimality condition for Problem~$(N)$ is given by,
\begin{equation}
\bzero \in \nu^0_f\,\partial E(\bx_f) + N_C(\bx_f)
\end{equation}
where, $\nu^0_f \ge 0$ is a Fritz John cost multiplier and $N_C(\bx_f)$ is the (limiting) normal cone to $C$ at $\bx_f$.  Now consider the following optimal control problem,
\begin{eqnarray}
&(M') \left\{
\begin{array}{lrl}
\textsf{Minimize }  & J[\bx(\cdot), \buf, t_f]
:=& E(\bx(t_f))  \\
\textsf{Subject to} & \dot\bx=& \bff'(\bx, \bu', t)\\
&(\bx(t_0), t_0) =& (\bx^0, t^0) \\
& \bx(t_f) \in & C \\
\end{array} \right.& \label{eq:prob-M'}
\end{eqnarray}
where, $\bu' \in \real{N_u}$ is a control variable, $\bff': \real{N_x} \times \real{N_u} \times \Real \to \real{N_x}$ is some given dynamics function, $t \in \Real$ is an independent ``time'' variable and $(\bx^0, t^0)$ is a given initial point in $\real{N_x} \times \Real$.  The terminal transversality condition for Problem~$(M')$ is given by,
\begin{equation}\label{eq:tvc-cone}
\blam'(t_f) \in \nu_0\,\partial E(\bx(t_f)) + N_C(\bx(t_f))
\end{equation}
where, $t_f$ is the final time, $\blam'(t_f) \in \real{N_x}$ is the final-time value of the adjoint covector and $\nu_0$ is the cost multiplier associated with \eqref{eq:prob-M'}.  Motivated by intellectual curiosity, a question posed in \cite{rossJCAM-1} was: Does an optimal control problem $(M') = (M)$ exist such that $\blam'(t_f) = \bzero$?  Needless to say, this question was answered in the affirmative for the case when $C$ is given by functional constraints,
\begin{equation}\label{eq:NLP-constraint}
C =  \set{\bx \in \real{N_x}:\ \be^L \le \be(\bx) \le \be^U }
\end{equation}
where, $\be: \bx \mapsto \real{N_e}, \ N_e \in \mathbb{N}$ is a given function, and $\be^L$ and $\be^U$ are the specified lower and upper bounds on the values of $\be$. Furthermore, the existence of Problem~$(M)$ was proved in \cite{rossJCAM-1} by direct construction.  No claim was staked on the uniqueness of such a problem.  In fact, the absence of uniqueness is utilized in this paper to devise another Problem~$(M)$ (in \Cref{sec:TMP}) that solves Problem~$(N)$.

It is apparent that the trajectory, $t \mapsto \bx(t)$, generated by Problem~$(M)$ is an ``algorithm'' for solving Problem~$(N)$, where $\bx(t_0) = \bx^0$ is the initial point or a guess to a solution for Problem~$(N)$. This observation implies that the traditional concept of an algorithm as a countable sequence generated by the point-to-set map,
\begin{equation}
\bx^0 \mapsto \set{\bx^0=\bx_0, \bx_1, \ldots, \bx_k, \bx_{k+1}, \ldots }
\end{equation}
be upgraded to its more primitive form:
\begin{equation}\label{eq:algorithm-mydef}
\bx^0 \mapsto \set{\bx^0 = \bx(t_0), [t_0, \infty) \ni t \mapsto \bx(t)}
\end{equation}
\begin{definition}[Algorithm Primitive]
Equation~\eqref{eq:algorithm-mydef} is an algorithm primitive for Problem~$(N)$.  A suitable discretization of \eqref{eq:algorithm-mydef} generates an algorithm given by,
\begin{equation}
\bx^0 \mapsto \set{\bx^0=\bx(t_0), \bx(t_1), \ldots, \bx(t_k), \bx(t_{k+1}), \ldots }
\end{equation}
\end{definition}
Suppose that an algorithm primitive is steerable by its tangent vector; then, we can write,
\begin{equation}\label{eq:ODE-single-int}
\dot\bx = \bu
\end{equation}
as a key equation that must constitute the vector field that defines Problem~$(M)$.  Although it was motivated by trajectory arguments, it is evident from a forward Euler discretization of \eqref{eq:ODE-single-int} that $\bu$ is, in fact, a continuous-time version of the search vector in optimization.  \emph{Note, however, that \eqref{eq:ODE-single-int} was not ``derived'' by considering the limit of a vanishing step size in optimization.}  In fact, it will be apparent later (in \Cref{sec:algorithms}) that there is a difference between an Eulerian and an optimization step-size.

Equation~\eqref{eq:ODE-single-int} was used in \cite{rossJCAM-1} to independently derive various algorithms such as the gradient and Newton's method.  Accelerated optimization algorithms appeared to be beyond the reach of the theory proposed in \cite{rossJCAM-1}; however, it was conjectured that such methods may be derivable by simply replacing \eqref{eq:ODE-single-int} by the double integrator model,
\begin{equation}\label{eq:ODE-double-int}
\ddot\bx = \bu
\end{equation}
The main contribution of this paper is in proving this conjecture.  A major consequence of this proof is a new approach to designing accelerated optimization algorithms.
\begin{remark}\label{rem:sobolev}
From an optimal control perspective, the difference between \eqref{eq:ODE-single-int} and \eqref{eq:ODE-double-int} seems quite trivial because the former implies $\bx(\cdot) \in W^{1,1}\left([t_0, t_f], \real{N_x}\right)$ while the latter indicates $\bx(\cdot) \in W^{2,1}\left([t_0, t_f], \real{N_x}\right)$\cite{vinter}.  Nonetheless, as will be apparent in the sections to follow, the ramifications of $\bx(\cdot)$ being an element of a smoother function space appear to have an outsized effect with regards to the problem of generating algorithms for solving Problem~$(N)$. From an optimization perspective, the differences between \eqref{eq:ODE-single-int} and \eqref{eq:ODE-double-int} is a little more nuanced: the search vector in \eqref{eq:ODE-single-int} steers the tangent vector (i.e., $\dot\bx$) whereas $\bu$ in \eqref{eq:ODE-double-int} steers the rate of change of the tangent vector.  Because the rate of change of the tangent vector implicitly incorporates prior information, the source of acceleration from the perspective of the algorithm primitive (i.e., $t \mapsto \bx(t)$) is in using this additional information to propel it forward. An interesting consequence of this observation is that algorithmic acceleration is indeed achieved by controlling acceleration (i.e., $\ddot\bx$).
\end{remark}
%

%

%

\section{A Transversality Mapping Principle}\label{sec:TMP}

With $C$ given by \eqref{eq:NLP-constraint}, the Lagrangian function for the nonlinear programming (NLP) Problem~$(N)$ may be written as,
\begin{equation}\label{eq:Lag-probN}
L(\nu^0_f, \bnu_f, \bx_f) := \nu^0_f E(\bx_f) +  \bnu_f\cdot\be(\bx_f)
\end{equation}
where, $(\nu^0_f, \bnu_f) \in \Real_+ \times \in \real{N_e}$ is the Fritz John multiplier pair, with $\bnu_f$ satisfying the complementarity condition, denoted by $\big(\bnu_f \dagger \be(\bx_f)\big)$, and given by,
\begin{equation}\label{eq:nu-comp}
\bnu_f \dagger \be(\bx_f) \quad \Leftrightarrow \quad \nu_i \left\{\begin{array}{ccrc}
               \le 0            & \text{if} & e_i(\bx_f) &= e_i^L \\
               =  0             & \text{if} & \qquad e_i^L < e_i(\bx_f) &< e_i^U \\
               \ge 0            & \text{if} & e_i(\bx_f) &= e_i^U \\
               unrestricted     & \text{if} & e_i^L &= e_i^U
             \end{array}
   \right.
\end{equation}
where, $i = 1, \ldots, N_x$. Together with \eqref{eq:nu-comp}, the first-order optimality condition for Problem~$(N)$ is given by,
\begin{equation}\label{eq:0=gradL}
\bzero = \partial_{\bx} L(\nu^0_f, \bnu_f, \bx_f)
\end{equation}
To construct Problem $(M)$, we follow \cite{rossJCAM-1} by ``sweeping back in time'' the data functions $E$ and $\be$ to define functions $t \mapsto y \in \Real$ and $t \mapsto \bs \in \real{N_e}$ according to,
\begin{subequations}\label{eq:idea-1}
\begin{align}
y(t) &:= E(\bx(t))\\
\bs(t) &:= \be(\bx(t))
\end{align}
\end{subequations}
Differentiating \eqref{eq:idea-1} with respect to time we get,
\begin{subequations}\label{eq:idea-2}
\begin{align}
\dot y  &= \big[\partial_{\bx} E(\bx)\big]\cdot \bv  \label{eq:cost-evolution}\\
\dot\bs &= \big[\partial_{\bx} \be(\bx)\big] \bv
\end{align}
\end{subequations}
where, we have set $\dot\bx := \bv$ as the ``velocity'' variable. Collecting all relevant equations, we construct the following time-free optimal control problem:
%
\begin{eqnarray}
&(M) \left\{
\begin{array}{lrl}
\textsf{Minimize }  & J[\bx(\cdot), \bv(\cdot), y(\cdot), \bs(\cdot), \buf, t_f]
:=& y(t_f)  \\
\textsf{Subject to} & \dot\bx=& \bv\\
&\dot\bv = & \bu \\
&\dot y=& \left[\partial_{\bx} E(\bx)\right]\cdot \bv\\
&\dot\bs =&  \left[\partial_{\bx}\be(\bx)\right] \bv\\
&\big(\bx(t_0), t_0\big) =& (\bx^0, t^0)\\
& \big(y(t_0), \bs(t_0)\big) = & \big(E(\bx^0), \be(\bx^0) \big)\\
& \bv(t_f) = & \bzero\\
&\be^L \le &\bs(t_f) \le  \be^U
\end{array} \right.& \label{eq:prob-M}
\end{eqnarray}
\vskip 6pt
\begin{remark}
Problem~$(N)$ is embedded in Problem~$(M)$. This follows from \eqref{eq:idea-1} and the imposition of the final-time constraint on $\bs(t)$ in \eqref{eq:prob-M}.  Furthermore, a solution to Problem~$(M)$ generates an algorithm primitive for Problem~$(N)$.
\end{remark}

The Pontryagin Hamiltonian\cite{vinter,ross-book} for Problem~$(M)$ is given by,
\begin{equation}\label{eq:H-def}
H(\blam_x, \blam_v, \lambda_y, \blam_s, \bx, \bv, y, \bs, \bu):= \blam_x\cdot\bv + \blam_v\cdot\bu + \lambda_y \left[\partial_{\bx} E(\bx)\right]\cdot \bv + \blam_s \cdot \left[\partial_{\bx}\be(\bx)\right] \bv
\end{equation}
where $\blam_x, \blam_v, \lambda_y$ and $\blam_s$ are the adjoint covectors corresponding to the dynamics associated with the variables $\bx, \bv, y $ and $\bs$ respectively.
\begin{lemma}\label{lemma:H-L-connection}
The Pontryagin Hamiltonian for Problem~$(M)$ and the instantaneous Lagrangian function associated with Problem~$(N)$ satisfy the condition,
\begin{equation}
H(\blam_x, \blam_v, \lambda_y, \blam_s, \bx, \bv, y, \bs, \bu):= \left[\blam_x + \partial_{\bx}L(\lambda_y, \blam_s, \bx)  \right]\cdot\bv + \blam_v \cdot \bu
\end{equation}
\end{lemma}
\begin{proof}
This follows directly from the defining equations given by \eqref{eq:Lag-probN} and \eqref{eq:H-def}.
\end{proof}
%
\begin{proposition}\label{prop:covec-evolve}
The adjoint arc $t \mapsto (\blam_x, \blam_v, \lambda_y, \blam_s)$ evolves according to,
\begin{subequations}
\begin{align}
\blam_x(t) & = - \partial_{\bx}L(\lambda_y(t), \blam_s(t), \bx(t)) + \bc_x\label{eq:covec-x-evolution} \\
\blam_v(t) & = -\bc_x (t - t_0) + \bc_v \label{eq:covec-v-evolution}\\
\lambda_y(t) & = c_y \label{eq:covec-y-evolution}\\
\blam_s(t) & = \bc_s \label{eq:covec-s-evolution}
\end{align}
\end{subequations}
where, $(\bc_x, \bc_v, c_y, \bc_s) \in \real{N_x} \times \real{N_x} \times \Real \times \real{N_e}$ is a constant.
\end{proposition}
\begin{proof}
The adjoint equations are given by,
\begin{subequations}
\begin{align}
\dot\blam_{x} &:=-\partial_{\bx}H = -\left[\partial^2_{\bx}L(\lambda_y, \blam_s, \bx)\right]\, \bv \label{eq:adj-x-NLP}\\
\dot\blam_v & := -\partial_{\bv}H =  -\blam_x - \partial_{\bx}L(\lambda_y, \blam_s, \bx) \label{eq:adj-v-NLP}\\
\dot\lambda_y & :=-\partial_y H = 0 \label{eq:adj-y-NLP} \\
\dot\blam_s & := -\partial_\bs H = \bzero \label{eq:adj-s-NLP}
\end{align}
\end{subequations}
Equations \eqref{eq:covec-y-evolution} and \eqref{eq:covec-s-evolution} follow directly from \eqref{eq:adj-y-NLP} and \eqref{eq:adj-s-NLP} respectively.

Substituting $\bv = \dot\bx$ in \eqref{eq:adj-x-NLP}, it follows that,
\begin{equation} \label{eq:lamxdot=Lagrangedot}
\dot\blam_{x} = -\frac{d}{dt}\Big[\partial_{\bx}L(\lambda_y, \blam_s, \bx)\Big]
+ \dot\lambda_y\, \partial_{\bx}E(\bx) + \sum_{i=1}^{N_e} \dot\lambda_{s_i} \partial_{\bx}e_i(\bx)
\end{equation}
Equation \eqref{eq:covec-x-evolution} follows from \eqref{eq:lamxdot=Lagrangedot}, \eqref{eq:adj-y-NLP} and \eqref{eq:adj-s-NLP}.

Substituting \eqref{eq:covec-x-evolution} in \eqref{eq:adj-v-NLP} we get
\begin{equation}
\dot\blam_v = - \bc_x
\end{equation}
from which \eqref{eq:covec-v-evolution} follows.
\end{proof}
\begin{theorem}\label{thm:H=0-infty}
\begin{enumerate}
\item[]
\item All extremals of Problem~$(M)$ are zero-Hamiltonian singular arcs.
\item All singular arcs of Problem~$(M)$ are of infinite order.
\end{enumerate}
\end{theorem}
\begin{proof}
From the Hamiltonian minimization condition we have the first-order condition,
\begin{equation}\label{eq:HMC4M=0}
\partial_\bu H(\blam_x, \blam_v, \lambda_y, \blam_s, \bx, \bv, y, \bs, \bu) = \bzero \Rightarrow \blam_v = \bzero \quad \forall\ t \in [t_0, t_f]
\end{equation}
Thus all extremals are singular.  From \cref{prop:covec-evolve} and \cref{eq:HMC4M=0} we get $\bc_x = \bzero$; hence, we have,
\begin{equation}\label{eq:lamx=-gradL}
\blam_x(t) = - \partial_{\bx}L(\lambda_y(t), \blam_s(t), \bx(t))
\end{equation}
The first part of the theorem now follows from \cref{lemma:H-L-connection}.

To prove the second part, differentiate $\partial_\bu H $ with respect to time:
\begin{equation}\label{eq:sing-arc-first-eq}
\frac{d}{dt}\partial_{\bu} H = \dot\blam_v =  -\blam_x - \partial_{\bx}L(\lambda_y, \blam_s, \bx)
\end{equation}
The second equality in \eqref{eq:sing-arc-first-eq} follows from \eqref{eq:adj-v-NLP}.  Differentiating \eqref{eq:sing-arc-first-eq} with respect to time we get,
\begin{align}
\frac{d^2}{dt^2}\partial_{\bu} H = \ddot\blam_v
&= - \dot\blam_x - \frac{d}{dt}\partial_{\bx}L(\lambda_y, \blam_s, \bx) \nonumber \\
& =- \dot\lambda_y\, \partial_{\bx}E(\bx) - \sum_{i=1}^{N_e} \dot\lambda_{s_i} \partial_{\bx}e_i(\bx)  \label{eq:sing-arc-second-eq}
\end{align}
where, the last equality follows from \eqref{eq:lamxdot=Lagrangedot}.  Substituting \eqref{eq:adj-y-NLP} and \eqref{eq:adj-s-NLP} in \eqref{eq:sing-arc-second-eq}, we get,
\begin{equation}
\frac{d^2}{dt^2}\partial_{\bu} H \equiv \bzero \quad \forall\ t \in [t_0, t_f]
\end{equation}
Hence, we have,
$$ \frac{d^k}{dt^k}\partial_{\bu} H = \bzero \quad\text{for\ } k = 0, 1 \ldots  $$
and no $k$ yields an expression for $\bu$.
\end{proof}

The endpoint Lagrangian\cite{ross-book} associated with the final-time conditions of Problem~$(M)$ may be written as,
\begin{equation}
\overline{E}(\nu_0, \bnu_v, \bnu_s, y(t_f), \bx(t_f), \bv(t_f), y(t_f), \bs(t_f)) := \nu_0 y(t_f) +\bnu_v \cdot \bv(t_f) + \bnu_s \cdot \bs(t_f)
\end{equation}
where, $\nu_0 \ge 0$ is the cost multiplier, $\bnu_v \in \real{N_x}$ and $\bnu_s$ satisfies the complementarity condition,
\begin{equation}\label{eq:nu-comp-proof}
\bnu_s \dagger \bs(t_f) \quad \Leftrightarrow \quad \nu_{s,i} \left\{\begin{array}{ccrc}
               \le 0            & \text{if} & s_i(t_f) &= e_i^L \\
               =  0             & \text{if} & \qquad e_i^L < s_i(t_f) &< e_i^U \\
               \ge 0            & \text{if} & s_i(t_f) &= e_i^U \\
               unrestricted     & \text{if} & e_i^L &= e_i^U
             \end{array}
   \right.
\end{equation}
Thus, the terminal transversality conditions for Problem~$(M)$ are given by,
\begin{subequations}
\begin{align}
\blam_{x}(t_f) &= \bzero \label{eq:tvc-x-NLP}\\
\blam_v(t_f) & = \bnu_v\\
\lambda_y(t_f) & = \nu_0 \ge 0 \label{eq:tvc-y-NLP}\\
\blam_s(t_f) & = \bnu_s \label{eq:tvc-s-NLP}
\end{align}
\end{subequations}
It is straightforward to show that the initial transversality generates the condition $\blam_v(t_0) = \bzero$; hence, $\bnu_v = \bzero$.
\begin{theorem}[Transversality Mapping Principle (TMP)]\label{thm:TMP}
The first-order necessary conditions for Problem~$(N)$ are imbedded in the terminal transversality conditions for Problem~$(M)$.
\end{theorem}
\begin{proof}
From \cref{thm:H=0-infty} (Cf.~\eqref{eq:lamx=-gradL}) and \eqref{eq:tvc-x-NLP} we get,
\begin{equation}\label{eq:tmp-proof-1}
\bzero = \partial_{\bx}L(\lambda_y(t_f), \blam_s(t_f), \bx(t_f))
\end{equation}
Substituting \eqref{eq:tvc-y-NLP} and \eqref{eq:tvc-s-NLP} in \eqref{eq:tmp-proof-1} we get the result (i.e., \eqref{eq:0=gradL} and \eqref{eq:nu-comp}) with the following mapping of the multipliers,
\begin{subequations}
\begin{align}
\nu^0_f &\longleftrightarrow \nu_0 = \lambda_y(t_f)\\
\bnu_f & \longleftrightarrow \bnu_s = \blam_s(t_f)
\end{align}
\end{subequations}
\end{proof}
\begin{remark}
\cref{thm:TMP} is an extension of the TMP presented in \cite{rossJCAM-1}.  Also, \cref{prop:covec-evolve} provides additional clarification and details that are absent in \cite{rossJCAM-1}.
\end{remark}

\section{New Principles for Accelerated Optimization}\label{sec:new-principles}
Because the extremals of Problem~$(M)$ are singular arcs of infinite order (Cf.~\Cref{thm:H=0-infty}), neither Pontryagin's Principle nor Krener's high order maximum principle\cite{krener-hmp} provide a computational mechanism for producing a singular optimal control.  Consequently, we need to develop new ideas for computation.

Collecting all the relevant primal-dual differential equations from \Cref{sec:TMP} together with their boundary conditions generates the following unconventional boundary value problem,
\begin{subequations}\label{eq:BVPfromTMP}
\begin{align}
\dot\bx & = \bv         & \dot\blam_{x} & = -\big[\partial^2_{\bx}L(\lambda_y, \blam_s, \bx)\big] \, \bv \\
\dot\bv & = \bu         & \dot\blam_v &  =  -\blam_x - \partial_{\bx}L(\lambda_y, \blam_s, \bx)  \\
\dot y & = \big[\partial_{\bx} E(\bx)\big] \cdot \bv    & \dot\lambda_y & = 0  \\
\dot\bs &=  \big[\partial_{\bx}\be(\bx)\big] \bv   & \dot\blam_s & = \bzero \\
\bx(t_0) & = \bx^0    &  \blam_x(t_f) & = \bzero \\
y(t_0) &= E(\bx^0)    &    \lambda_y(t_f) & \ge 0\\
\bs(t_0) &= \be(\bx^0) &    \blam_s(t_f) &\dagger \bs(t_f) \\
\blam_v(t_0) & = \bzero  & \bv(t_f) & = \bzero   \\
              &           & \be^L &\le \bs(t_f) \le  \be^U
\end{align}
\end{subequations}
Any infinite-order singular control trajectory $\buf$ that solves \eqref{eq:BVPfromTMP} also solves Problem~$(M)$. Consequently, such a solution generates an algorithm primitive that solves Problem~$(N)$.  To produce such an algorithm primitive, we follow and extend the ideas proposed in \cite{rossJCAM-1} by using the sweeping principle to inject dual control variables.  That is, as in \cite{rossJCAM-1}, we replace the equation $\dot\blam_s = \bzero $ by introducing a control variable $\bmu$ that steers $\blam_s(t)$:
\begin{equation}\label{eq:inject-1}
\dot\blam_s = \bmu
\end{equation}
Similarly, we set,
\begin{equation}\label{eq:inject-2}
\dot\lambda_y = \omega
\end{equation}
In the unaccelerated version of this theory\cite{rossJCAM-1}, it was important modify the adjoint equation (corresponding to $\bx$) to maintain a zero-Hamiltonian singular trajectory (Cf.~\Cref{thm:H=0-infty}).  Adopting the same idea, we modify the adjoint equation according to,
\begin{equation}\label{eq:lam_x-mod}
-\dot\blam_{x} = \big[\partial^2_{\bx}L(\lambda_y, \blam_s, \bx)\big] \, \bv + \partial_{\bx}L(\omega, \bmu, \bx)
\end{equation}
Equation~\cref{eq:lam_x-mod} is simply the time derivative of \eqref{eq:lamx=-gradL}.  Consequently, \cref{eq:lam_x-mod} also ensures that $\dot\blam_v = \bzero \, \forall\ t \in [t_0, t_f]$ (Cf.~\eqref{eq:adj-v-NLP}); hence, $\blam_v$ can be safely eliminated in generating a singular solution to \eqref{eq:BVPfromTMP}.  Thus, the problem of generating a candidate infinite-order singular arc to Problem~$(M)$ reduces to a controllability-type problem associated with the following auxiliary primal-dual system,
\begin{equation}\label{eq:A4auxSystem}
(A) \left\{
\begin{aligned}
\dot\blam_{x}  &= -\big[\partial^2_{\bx}L(\lambda_y, \blam_s, \bx)\big] \, \bv - \partial_{\bx}L(\omega, \bmu, \bx)\\
\dot\lambda_y &=\omega\\
\dot\blam_s & = \bmu\\
\dot\bv & = \bu\\
\dot\bs &=  \big[\partial_{\bx}\be(\bx)\big] \bv
\end{aligned}
\right.
\end{equation}
The final-time conditions for $(A)$ are extracted from \eqref{eq:BVPfromTMP} and can be specified in terms of the target set, $T$ given by,
\begin{multline}\label{eq:target4A}
T := \left\{\blam_x(t_f), \lambda_y(t_f), \blam_s(t_f), \bv(t_f), \bs(t_f)\mid \right.\\
\left.\blam_x(t_f) = \bzero,\ \lambda_y(t_f) \ge 0, \ \blam_s(t_f) \dagger \bs(t_f), \
\bv(t_f) = \bzero, \ \be^L \le \bs(t_f) \le  \be^U   \right\}
\end{multline}

In the discussions to follow, it will be convenient to view the dynamical system $(A)$ in terms of the sum of two vector fields:
\begin{equation}\label{eq:f=sumOf2}
\dot\bz = \bff(\lambda_y, \blam_s, \bv, \bx, \bzeta) :=
\underbrace{
\begin{bmatrix}
-\big[\partial^2_{\bx}L(\lambda_y, \blam_s, \bx)\big] \, \bv \\
0\\
\bzero\\
\bzero\\
\big[\partial_{\bx}\be(\bx)\big] \bv
\end{bmatrix}
}_{\bff_0}
+
\underbrace{
\begin{bmatrix}
- \partial_{\bx}L(\omega, \bmu, \bx)\\
\omega\\
\bmu\\
\bu\\
\bzero
\end{bmatrix}
}_{\bff_1}
\end{equation}
where,
\begin{subequations}
\begin{align}
\bz:= & (\blam_x, \lambda_y, \blam_s, \bv, \bs) & \bzeta:= & (\bu, \bmu, \omega) \\
\bff_0 \equiv & \bff_0(\lambda_y, \blam_s, \bv, \bx) & \bff_1 \equiv &\bff_1(\bx, \bzeta)
\end{align}
\end{subequations}
In control theory, $\bff_0$ is known as the drift vector field, whose presence (or absence) impact the production of solutions to the $(A)$-$(T)$ system. In the unaccelerated version of this theory, there is no drift vector field\cite{rossJCAM-1}; hence, an extension of the ideas to accelerated optimization requires an explicit consideration of $\bff_0$.

The main problem of interest with respect to generating an algorithm primitive for solving Problem~$(N)$ can now be framed as finding the control function $t \mapsto \bzeta $ that drives a given point, $T^c \ni \bz^0 = \bz(t_0)$, to some point $\bz(t_f) \in T$, where, $T^c$ is the complement of $T$.  To formalize the statement of this problem, we adopt Clarke's notion of guidability\cite{clarkeLyap,clarkeEncy}:
\begin{definition}[Guidability]\label{def:guidability}
A point $\bz^0 \in T^c$ is guidable to $T$ if there is a trajectory $[t_0, t_f] \to \bz(t)$ satisfying $\bz(t_0) = \bz^0$ and $\bz(t_f) \in T$.
\end{definition}
\begin{definition}[Global Guidability]\label{def:guidability-global}
A point $\bz^0 \in T^c$ is globally guidable to $T$ if every point $\bz^0 \in T^c$ is guidable to $T$.
\end{definition}
\begin{definition}[Asymptotic Guidability]\label{def:guidability-asymp}
A point $\bz^0$ is asymptotically guidable to $T$ if it is guidable with $t_f \rightarrow \infty$.
\end{definition}
%

It is apparent that the notion of guidability is weaker than stability.  Furthermore, it is clear that guidability is quite sufficient in terms of producing an algorithm primitive to solve Problem~$(N)$.

To design algorithm primitives for Problem~$(N)$, we simply need to find guidable trajectories for the $(A)$-$(T)$ pair.  A standard workhorse in control theory for solving such a problem is a control Lyapunov function (CLF)\cite{clarkeLyap,sontag-book}.  Following \cite{motta-CLFs}, we define a CLF for the $(A)$-$(T)$ system as a positive definite function, $V:\overline{T^c} \to \Real$, such that for each point in $T^c$, there exists a value of $\bff$ that points in a direction along which $V$ is strictly decreasing.  Let $\pounds_f V$ be the Lie derivative of $V$ along the vector field $\bff$.  Then the strict decreasing condition can be expressed as,
\begin{equation}\label{eq:clf-theory-1}
\pounds_f V := \Big\langle\partial V(\bz),\  \bff( \lambda_y, \blam_s, \bv, \bx, \bzeta) \Big\rangle < 0
\end{equation}
for some choice of $\bzeta$.  Because it is possible for $\pounds_{f_1} V $ to vanish for all choices of $\bzeta$ (see \eqref{eq:f=sumOf2}) when $\bz \in T^c$, a satisfaction of \eqref{eq:clf-theory-1} requires the condition,
\begin{equation}\label{eq:clf-theory-2}
\pounds_{f_0} V := \Big\langle\partial V(\bz),\  \bff_0( \lambda_y, \blam_s, \bv, \bx) \Big\rangle < 0 \quad \text{if } \pounds_{f_1} V = 0
\end{equation}
whenever $\bz \not\in T$.

As a means to get the best instantaneous solution, suppose we choose controls such that
\begin{equation}\label{eq:minP-not}
\bzeta = \arg\min_{\bzeta} \pounds_f V
\end{equation}
One problem with \eqref{eq:minP-not} is that $\pounds_f V$ is an affine function of $\bzeta$ and the control is unbounded. Hence, to use \eqref{eq:minP-not} in a meaningful manner, it is necessary to constrain $\bzeta$ to some compact set $\U$.  This notion is similar to that of a trust region in optimization; however, as shown in \cite{rossJCAM-1}, a proper choice for $\U$ also generates new insights on the selection of a metric space for an optimization algorithm.  Hence, we frame the idea implicit in \eqref{eq:minP-not} in terms of the following minimum principle:
\begin{eqnarray}\label{eq:MinP4search}
&\quad (P) \left\{
\begin{array} {lll}
\displaystyle\mathop\textsf{Minimize }_{\bzeta}  & \pounds_f V := \Big\langle\partial V(\bz),  \bff( \lambda_y, \blam_s, \bv, \bx, \bzeta) \Big\rangle\\
\textsf{Subject to} & \bzeta \in  \U (\bz, \bx, y, t)
\end{array} \right.&
\end{eqnarray}
where, $\U (\bz, \bx, y, t)$ is any given compact set that may vary with respect to the tuple $(\bz, \bx, y, t)$; i.e.,  $\U : (\bz, \bx, y, t) \rightrightarrows \real{N_x} \times \real{N_e} \times \Real$.  Equation~\eqref{eq:MinP4search} is a direct extension of the minimum principle posed in \cite{rossJCAM-1}.  The caveat in applying \eqref{eq:MinP4search} is an assurance of \eqref{eq:clf-theory-2}.

In exploring a different method to manage the drift vector field, we exchange the cost function and constraint condition in \eqref{eq:MinP4search} to formulate an alternative minimum principle that holds the potential to provide additional insights in formulating optimal algorithm primitives.  To facilitate this development, we select $\rho:(\bz, \bx, y, t) \mapsto \Real_+$ to be some function such that $-\rho$ specifies a rate of descent for $\pounds_f V$.  That is, we replace \eqref{eq:clf-theory-1} by the constraint,
\begin{equation}\label{eq:preMinP*}
\exists\, \bzeta \text{ s.t. } \pounds_f V := \Big\langle\partial V(\bz),\  \bff( \lambda_y, \blam_s, \bv, \bx, \bzeta) \Big\rangle
\le -\rho(\bz, \bx, y, t)
\end{equation}
Let  $D: (\bzeta, \bz, \bx, y, t) \mapsto \Real$ be an appropriate objective function. Then, an alternative minimum principle may be posed as:
\begin{eqnarray}\label{eq:MinP*4search}
&(P^*) \left\{
\begin{array} {lll}
\displaystyle\mathop\textsf{Minimize }_{\bzeta}  &D (\bzeta, \bz, \bx, y, t)  \\
\textsf{Subject to} & \pounds_f V +\rho(\bz, \bx, y, t) \le 0
\end{array} \right.&
\end{eqnarray}
An apparently obvious choice for $\rho$ in \eqref{eq:MinP*4search} is $V$ itself because it would imply that the resulting Lyapunov function would decrease at least exponentially. As fast as an exponential might be, it turns out a better choice for $\rho$ may be possible if we view the minimum principles $(P)$ and $(P^*)$ as merely computational techniques to solve the CLF inequality\cite{clarkeLyap},
\begin{equation}\label{eq:CLF-inequality}
\min_{\scriptsize{\bzeta} \in \U} \pounds_f V +\rho(\bz, \bx, y, t) \le 0
\end{equation}
As is well documented\cite{clarkeLyap, motta-CLFs,CLFtoHJB2020}, what is most interesting about \eqref{eq:CLF-inequality} is that it can be rewritten as a Hamilton-Jacobi-Bellman (HJB) inequality,
\begin{equation}\label{eq:HJB-inequality}
\min_{\scriptsize{\bzeta} \in \U} H^P(\partial V, \bz, \bx, \bzeta) +\rho(\bz, \bx, y, t) \le 0
\end{equation}
where, $H^P(\blam^A, \bz, \bx, \bzeta) := \langle\blam^A, \bff(\lambda_y, \blam_s, \bv, \bx, \bzeta)\rangle$ may be viewed as the Pontryagin Hamiltonian for System~$(A)$. Evidently, even a minimum-time solution can be produced if $V$ is chosen as the time-to-go function\cite{clarkeLyap}.
Because such ``optimal functions'' are unknown, a more tractable approach to selecting $\rho$ is provided by the following theorem due to Bhat and Bernstein\cite{bhat-2000}:
\begin{theorem}[Bhat-Bernstein]\label{thm:BB}
Let $\rho$ be given by,
\begin{equation}
\rho(\bz) := r \big(V(\bz)\big)^{1-m}
\end{equation}
where $r > 0$ and $m \in (0,1)$. Then, the the time interval for a guidable trajectory $[t_0, t_f] \mapsto \bz$  is bounded by,
\begin{equation}
(t_f-t_0) \le \frac{\big(V(\bz^0)\big)^m}{r\,m}
\end{equation}
\end{theorem}
%
\begin{remark}
It is apparent that the ``left'' limiting case of $m \to 0$ in \cref{thm:BB} corresponds to the case of asymptotic guidability while the ``right'' limiting case of $m \to 1$ may be viewed as a solution to a  minimum-time problem provided $V$ is chosen as the time-to-go function\cite{clarkeLyap}.
\end{remark}
\begin{remark}
Based on the connections between the HJB equations and a CLF as a computational method for selecting $\bzeta$, the minimum principles $(P)$ and $(P^*)$ may be viewed as Pontryagin-type conditions for an optimal control of System~$(A)$.
\end{remark}

The minimum principles $(P)$ and $(P^*)$ are technically not new.  They have been widely used in control theory for generating feedback controls\cite{clarkeLyap,sontag-book,motta-CLFs,freeman-acc}.  What makes them new in \eqref{eq:MinP4search} and \eqref{eq:MinP*4search} is their specific use for the $(A,T)$ pair, and consequently, in designing ordinary differential equations (ODEs) that generate algorithm primitives (cf.~\eqref{eq:algorithm-mydef}).  Furthermore, recall that the $(A,T)$ system was derived from the necessary conditions of Problem~$(M)$ with the TMP providing the critical link (Cf.~\cref{thm:TMP}) between Problems~$(M)$ and $(N)$. These ideas are in sharp contrast to earlier works\cite{NLP2ODE-1980, NLP2ODE-1981, NLP2ODE-1989, NLP2ODE-1994, NLP2ODE-2006, NLP2ODE-2007, NLP2ODE-2017} that have sought to solve NLPs using differential equations.  Consequently, the differential equations proposed in these prior works are not only different from \eqref{eq:A4auxSystem} but also that \emph{we use \eqref{eq:A4auxSystem} as generators of ODEs}.  Furthermore, the CLFs used in \eqref{eq:MinP4search} and \eqref{eq:MinP*4search} are generic; hence, different choices of $V$ can lead to different ODEs, which, in turn generate different algorithm primitives. Finally, note also that the focus of the current ideas is primarily on accelerated optimization.

In the absence of additional analysis, it might appear that we have come to full circle; i.e., in the quest for solving NLPs via optimal control theory, we have generated Problems~$(P)$ and $(P^*)$ that appear to be NLPs themselves. As a result, the proposed theory would only be meaningful if (a) it lead to some new insights on solving Problem~$(N)$ and/or (b) Problems~$(P)$ and $(P^*)$ were simpler than $(N)$. Because the unaccelerated version of this theory\cite{rossJCAM-1} did indeed generate new insights, the same can be expected in pursuing this idea further. This is shown in \Cref{sec:derive-HB+}.  In addition, because System~$(A)$ is affine in the control variable, Problems~$(P)$ and $(P^*)$ can indeed be rendered simpler than $(N)$.  In this context we briefly note that the structure of the vector field $\bff$ can be further altered quite easily through the process of adding more integrators.  For example, analogous to \eqref{eq:ODE-double-int}, we can replace $\dot\blam_s = \bmu $ and $\dot\lambda_y = \omega$ by,
\begin{align}
\dot\blam_s = \btheta_s, \quad \dot\btheta_s = \bomega_s \\
\dot\lambda_y = \theta_y, \quad \dot\theta_y = \bomega_y
\end{align}
to generate a new $(A, T)$ pair:
\begin{eqnarray}
\qquad (A') \left\{
\begin{aligned}
\dot\blam_{x}  &= -\big[\partial^2_{\bx}L(\lambda_y, \blam_s, \bx)\big] \, \bv - \partial_{\bx}L(\theta_y, \btheta_s, \bx)\\
\dot\lambda_y &=\theta_y\\
\dot\theta_y & = \omega_y\\
\dot\blam_s & = \btheta_s\\
\dot\btheta_s & = \bomega_s\\
\dot\bv & = \bu\\
\dot\bs &=  \big[\partial_{\bx}\be(\bx)\big] \bv
\end{aligned}
\right.
& (T') \left\{
\begin{aligned}
\blam_x(t_f) & = \bzero\\
\lambda_y(t_f) &\ge 0 \\
\theta_y(t_f) & = 0\\
\btheta_s(t_f) & = \bzero\\
\bv(t_f) & = \bzero\\
\be^L \le &\bs(t_f) \le  \be^U  \\
\blam_s(t_f) &\dagger \bs(t_f)
\end{aligned}
\right.
\end{eqnarray}
%
As noted earlier (Cf.~\cref{rem:sobolev}), the addition of integrators seems to have a profound effect on the production of accelerated algorithms.

\section{Generation of Accelerated Algorithm Primitives Illustrated}\label{sec:derive-HB+}
To illustrate some specific features of the general theory presented in \Cref{sec:TMP} and \Cref{sec:new-principles}, consider the unconstrained optimization problem,
\begin{equation}
(S) \left\{\displaystyle\mathop\text{Minimize }_{\bx_f \in \real{N_x}}   E(\bx_f) \right.
\end{equation}
Producing accelerated algorithms for such problems have generated increased attention in recent years \cite{lessard-2016,su-2016,wibisono-2016} due to their immediate applicability to machine learning.

\subsection{Development of the Auxiliary System}

From \Cref{sec:TMP}, it follows that the optimal control problem that solves Problem~$(S)$ is given by:
\begin{eqnarray}
&(R) \left\{
\begin{array}{lrl}
\textsf{Minimize }  & J[y(\cdot), \bx(\cdot), \bv(\cdot), \buf, t_f]
:=& y_f  \\
\textsf{Subject to} & \dot\bx=& \bv\\
&\dot\bv = & \bu \\
&\dot y=& \left[\partial_{\bx} E(\bx)\right]\cdot \bv\\
&(\bx(t_0), t_0) =& (\bx^0, t^0) \\
& y(t_0) = & E(\bx^0) \\
&\bv(t_f) = & \bzero
\end{array} \right.& \label{eq:prob-R}
\end{eqnarray}
\begin{remark}
The unaccelerated version of Problem~$(R)$ (i.e., one without the velocity variable, $\bv$) was first formulated by Goh\cite{Goh-1997}; however, because the problem is singular (cf.~\Cref{thm:H=0-infty}), Goh et al\cite{Goh-2021} advanced an alternative theory based on bang-bang controls by adding control constraints to (the unaccelerated version of) Problem~$(R)$.
\end{remark}
\begin{proposition}\label{prop:normality}
Problem~$(R)$ has no abnormal extremals.
\end{proposition}
\begin{proof}
This proof is straightforward; hence, it is omitted.
\end{proof}

It is straightforward to show that \cref{eq:BVPfromTMP} reduces to,
\begin{subequations}\label{eq:BVP4R}
\begin{align}
\dot\bx & = \bv         & \dot\blam_{x} & = -\lambda_y\,\partial^2_{\bx}E(\bx) \, \bv \\
\dot\bv & = \bu         & \dot\blam_v & =  -\blam_x - \lambda_y\, \partial_{\bx} E(\bx)\\
\dot y & = \big[\partial_{\bx} E(\bx)\big]^T \bv    & \dot\lambda_y & = 0\\
\bx(t^0) & = \bx^0  & \bv(t_f) & = \bzero\\
y(t^0) &= E(\bx^0)    & \blam_x(t_f) & = \bzero\\
\blam_v(t^0) & =  \bzero  & \lambda_y(t_f) & > 0
\end{align}
\end{subequations}
It thus follows that the auxiliary primal-dual dynamical system is given by,
\begin{equation}\label{eq:aux-dynamics-R}
(A_R) \left \{
\begin{aligned}
\dot\blam_{x} & = -\partial^2_{\bx}E(\bx) \, \bv   \\
\dot\bv & = \bu
\end{aligned}
\right.
\end{equation}
where, we have scaled the adjoint covector by the constant, $\lambda_y > 0$ (Cf.~\cref{prop:normality}).  The target final-time condition for $(A_R)$ is given by,
\begin{equation}\label{eq:aux-target-R}
(T_R) \left \{
\begin{aligned}
\blam_x(t_f) & = \bzero\\
\bv(t_f) & = \bzero
\end{aligned}
\right.
\end{equation}
\subsection{Application of the Minimum Principles}
Following \eqref{eq:f=sumOf2} we write $\bff$ for $(A_R)$ as,
\begin{equation}
\bff(\bx, \bv, \bu) := \underbrace{\left[
                        \begin{array}{c}
                          -\partial^2_{\bx}E(\bx) \, \bv\\
                           \bzero\\
                        \end{array}
                      \right] }_{\bff_0}
                    + \underbrace{\left[
                        \begin{array}{c}
                          \bzero\\
                           \bu\\
                        \end{array}
                      \right] }_{\bff_1}
\end{equation}
Furthermore, if $V:(\blam_x, \bv) \mapsto \Real$ is a CLF, then we must have,
\begin{equation}\label{eq:CLF-condition-2}
\pounds_f V =  -\big\langle\partial_{\boldsymbol\lambda_x}V(\blam_x, \bv),\ \partial^2_{\bx}E(\bx) \, \bv \big\rangle + \big\langle\partial_{\bv}V(\blam_x, \bv),\ \bu \big\rangle< 0
\end{equation}
for some choice of $\bu$  whenever $(\blam_x, \bv) \ne (\bzero, \bzero)$.  In addition, \eqref{eq:clf-theory-2} simplifies to,
\begin{equation}\label{eq:drift-condition}
\big\langle \partial_{\boldsymbol\lambda_x}V(\blam_x, \bv),\ \partial^2_{\bx}E(\bx) \, \bv \big\rangle > 0 \quad\text{if }\quad \partial_{\bv}V(\blam_x, \bv) = \bzero \ \text{and } (\blam_x, \bv) \ne (\bzero, \bzero)
\end{equation}
Furthermore, we set $\bu = \bzero$ if $\partial_{\bv} V = \bzero$. This last statement implies that the dynamical system $(A_R)$ will continue to evolve as a result of $\bv \ne \bzero$.

Let $\U (\bx, \blam_x, \bv, t)$ be a compact set that may vary with respect to the tuple $(\bx, \blam_x, \bv, t)$;  then, \eqref{eq:MinP4search} may be formulated as,
\begin{eqnarray}
&(P_S) \left\{
\begin{array} {lll}
\displaystyle\mathop\textsf{Minimize }_{\bu}  & \pounds_f V := \big\langle\partial V(\blam_x, \bv),\ \bff(\bx, \bv, \bu) \big\rangle \\
\textsf{Subject to} & \bu \in  \U (\bx, \blam_x, \bv, t)
\end{array} \right.& \label{prob:Min-P4S}
\end{eqnarray}
To formulate Problem~$(P^*_S)$ that is analogous to \eqref{eq:MinP*4search}, we select a function  $D: (\bu, \bx, \blam_x, \bv) \mapsto \Real$ to be an appropriate objective function.  Then, an application of \eqref{eq:MinP*4search} reduces to,
\begin{eqnarray}
&(P^*_S) \left\{
\begin{array} {lll}
\displaystyle\mathop\textsf{Minimize }_{\bu}  &D (\bu, \bx, \blam_x, \bv, t) \\
\textsf{Subject to} & \pounds_f V +\rho(\blam_x, \bv, \bx, t) \le 0
\end{array} \right.& \label{prob:Min-P*}
\end{eqnarray}
The generation of accelerated algorithm primitives is now reduced to designing $V$ and $\U$ in $(P_S)$ or $D, V$ and $\rho$ in $(P^*_S)$.

\subsection{Optimal Control for Some Accelerated Algorithm Primitives}
Let $\bW : (\bx, \blam_x, \bv, t) \mapsto \mathbb{S}^{N_x}_{++}$ be a symmetric positive definite matrix function that metricizes the space $\U$. Following \cite{rossJCAM-1}, we consider
\begin{equation}\label{eq:U-metric-trust}
\U(\bx, \blam_x, \bv, t) := \set{\bu:\  \bu^T\bW(\bx, \blam_x, \bv, t)\bu \le \Delta(\bx, \blam_x, \bv, t) }
\end{equation}
where $\Delta: (\bx, \blam_x, \bv, t) \mapsto  \Real_{++}$.  Note that $\Delta$ is similar to, but is not, the familiar trust region in optimization. Under these conditions, a solution to \eqref{prob:Min-P*} is given explicitly by,
\begin{equation}\label{eq:u-fromMinP}
\bu = \left\{
        \begin{array}{ll}
          -\sigma[@t]\, \bW^{-1}[@t] \, \partial_{\bv} V(\blam_x, \bv) & \hbox{ if\ } \partial_{\bv}V(\blam_x, \bv) \ne \bzero\\
          \bzero & \hbox{ if\ } \partial_{\bv}V(\blam_x, \bv) = \bzero
        \end{array}
      \right.
\end{equation}
where,
\begin{equation}\label{eq:sigma-def}
\sigma[@t] := + \sqrt{
\frac{\Delta(\bx, \blam_x, \bv, t)}{\big[\partial_{\bv}V(\blam_x, \bv)\big]^T \bW^{-1}[@t]\big[\partial_{\bv}V(\blam_x, \bv)\big]}
}
\end{equation}
and $\bW[@t]\equiv \bW(\bx, \blam_x, \bv, t)$.

To illustrate an application of Minimum Principle $(P^*)$, we select
\begin{equation}\label{eq:U-metric-trust*}
D (\bu, \bx, \blam_x, \bv, t) =\frac{1}{2} \big(\bu^T\bW(\bx, \blam_x, \bv, t)\bu\big)
\end{equation}
Solving the resulting problem, we get
\begin{equation}\label{eq:u-fromMinP*}
\bu = \left\{
        \begin{array}{ll}
          -\sigma^*[@t]\, \bW^{-1}[@t] \, \partial_{\bv} V(\blam_x, \bv) & \hbox{ if\ } \partial_{\bv}V(\blam_x, \bv) \ne \bzero\\
          \bzero & \hbox{ if\ } \partial_{\bv}V(\blam_x, \bv) = \bzero
        \end{array}
      \right.
\end{equation}
where,
\begin{equation}\label{eq:sigma*-def}
\sigma^*[@t] := \left\{
        \begin{array}{ll}
          \frac{\xi({\boldsymbol\lambda}_x, \bv, \bx, t)}{\big[\partial_{\bv}V(\blam_x, \bv)\big]^T \bW^{-1}[@t]\big[\partial_{\bv}V({\boldsymbol\lambda}_x, \bv)\big]} & \hbox{ if\ } \xi({\boldsymbol\lambda}_x, \bv, \bx) \ge 0\\
          0 & \hbox{ if\ } \xi(\blam_x, \bv, \bx) < 0
        \end{array}
      \right.
\end{equation}
and
\begin{equation}
\xi(\blam_x, \bv, \bx, t):= \rho(\blam_x, \bv, \bx, t)  -\big\langle\partial_{\boldsymbol\lambda_x}V(\blam_x, \bv),\ \partial^2_{\bx}E(\bx) \, \bv \big\rangle
\end{equation}
Comparing \eqref{eq:u-fromMinP} and \eqref{eq:u-fromMinP*} it follows that for the choice of $\U$ and $D$ given by \eqref{eq:U-metric-trust} and \eqref{eq:U-metric-trust*} respectively, both minimum principles ($P$ and $P^*$) generate the same functional form for $\bu$ but with different interpretations for the control ``gains'' given by $\sigma$ and $\sigma^*$.


\subsection{Generation of ODEs For Some Accelerated Optimization Algorithms}

\begin{proposition}
Let,
\begin{equation}\label{eq:V-quad-pd}
V(\blam_x, \bv) =  (a/2)\blam_x^T\blam_x + (b/2)\bv^T\bv+ c\blam_x^T\bv
\end{equation}
where, $a > 0, \quad b > 0$ and $ c < 0$ are real numbers such that $ab - c^2 > 0$.  Then, if $E$ is a strictly convex function, $V(\blam_x, \bv)$ is a CLF for the $(A_R)$-$(T_R)$ pair.
\end{proposition}
\begin{proof}
The conditions $a >0, b >0$ and $ab-c^2 > 0$ ensure that $V$ is positive definite.  The Lie derivative of $V$ along $\bff$ is given by,
\begin{equation}
\pounds_f V = \big\langle a \blam_x + c \bv,\  -\partial^2_{\bx} E(\bx) \bv  \big\rangle  + \big\langle c \blam_x + b \bv,\ \bu  \big\rangle
\end{equation}

If $c \blam_x + b \bv \ne \bzero$, then choosing $\bu$ according to \eqref{eq:u-fromMinP*} ensures that $\pounds_f V < 0$ for any choice of $\rho(\blam_x, \bv, \bx) > 0$.

If $c \blam_x + b \bv = \bzero$, then $\pounds_{f_1} V = 0$ for all choices of $\bu$. In this case, $\blam_x = -(b/c) \bv$; hence, we have
\begin{align}
\pounds_{f_0} V  & = \big\langle a\blam_x + c \bv,\ -\partial^2_{\bx}E(\bx) \, \bv \big\rangle \nonumber \\
                &= \left(\frac{ab - c^2}{c}\right)\bv^T \partial^2_{\bx}E(\bx) \, \bv \ < 0 \quad \text{if\ } (\blam_x, \bv) \ne (\bzero, \bzero) \label{eq:clf-theory-2-proof4quad}
\end{align}
where, the inequality in \eqref{eq:clf-theory-2-proof4quad} follows from $c < 0$ and $E$ strictly convex.  Hence, $V$ satisfies \eqref{eq:clf-theory-2}.
\end{proof}
\begin{corollary}
Polyak's equation\cite{polyak64} for the heavy ball method can be generated from the minimum principles $(P_S)$ or $(P^*_S)$ using a Euclidean metric for $\bW$ and the quadratic CLF given by \eqref{eq:V-quad-pd}.
\end{corollary}
\begin{proof}
Let $\sigma^q$ denote $\sigma$ or $\sigma^*$ given by \eqref{eq:sigma-def} and \eqref{eq:sigma*-def} respectively. Let,
\begin{subequations}\label{eq:gamma-def}
\begin{align}
\gamma^a[@t] &:= -c\,\sigma^q[@t] \ge 0\\
\gamma^b[@t] &:= b\,\sigma^q[@t] \ge 0
\end{align}
\end{subequations}
Using \eqref{eq:V-quad-pd}, the expression for $\bu$ given by either \eqref{eq:u-fromMinP} or \eqref{eq:u-fromMinP*} can be written universally as,
\begin{equation} \label{eq:uFamilyNo1}
\bu = - \bW^{-1}[@t]\big(\gamma^a[@t]\,\partial_{\bx}E(\bx) + \gamma^b[@t]\, \bv\big)
\end{equation}
where, we have used the integral of motion $\blam_x = -\partial_{\bx} E(\bx)$ in accordance with \eqref{eq:lamx=-gradL}.  Substituting \eqref{eq:uFamilyNo1} in  \eqref{eq:ODE-double-int} we get,
\begin{equation}\label{eq:ODE-HB++}
\bW[@t]\, \ddot\bx + \gamma^a[@t]\,\partial_{\bx}E(\bx) + \gamma^b[@t]\, \bv = \bzero
\end{equation}
Polyak's equation is given by\cite{polyak64},
\begin{equation}\label{eq:ODE-HB}
\ddot\bx + a_1(t)\, \dot\bx + a_2(t)\,\partial_{\bx}E(\bx) = \bzero
\end{equation}
where $a_1(t) > 0 $ and $a_2(t) > 0$ are time-varying scalar parameters. Equation~\eqref{eq:ODE-HB} thus follows from \eqref{eq:ODE-HB++} with $\bW$ set to the identity matrix.
\end{proof}
\begin{remark}\label{rem:HB++}
Polyak ``derived'' \eqref{eq:ODE-HB} based on physical considerations of the motion of ``a small heavy sphere''\cite{polyak64}. A discrete analog of \eqref{eq:ODE-HB} generates his momentum method. In \cite{polyak2017}, Polyak et al argue that \eqref{eq:ODE-HB} also generates Nesterov's accelerated gradient method\cite{nesterov83} if $a_1(t)$ is set to $3/t$.  This specific choice of $a_1(t)$ is based on the results of Su et al\cite{su-2016}.
\end{remark}

From \Cref{rem:HB++} it follows that \eqref{eq:ODE-HB++} can generate both Polyak's momentum method and Nesterov's accelerated gradient method.  Evidently, alternative accelerated optimization algorithms are possible by various selection of the parameters in \eqref{eq:ODE-HB++}.

\section{A New Approach to Generating Algorithms}\label{sec:algorithms}
The results of \Cref{sec:derive-HB+} demonstrate that the minimum principles $(P)$ and $(P^*)$ can successfully generate ODEs that govern the flow of accelerated algorithm primitives.  It thus seems reasonable to suggest that algorithms can be produced by simply discretizing the resulting ODEs.  We depart from this perspective for a variety of reasons, some of which are implied in \Cref{rem:HB++}.  To clarify the need for a new approach to generating algorithms, consider a discretization of \eqref{eq:ODE-HB++} with $\bW$ set to the identity matrix.  From elementary numerical methods, it is straightforward to produce the following algorithm:
\begin{equation}\label{eq:HB-gain-step}
\bx_{k+1} = \bx_k - \underbrace{\left(h^2_k \gamma^a_k \right)}_{\alpha_k} \partial_{\bx}E(\bx_k) + \underbrace{\left(1 - h_k\gamma^b_k\right)}_{\beta_k} (\bx_k - \bx_{k-1})
\end{equation}
Equation~\eqref{eq:HB-gain-step} indicates the connections between a discretization step, $h_k$, associated with \eqref{eq:ODE-HB++}, the step length, $\alpha_k$, in optimization, the momentum parameter, $\beta_k$, associated with the heavy-ball method and the discretized controller gains $\gamma^a_k$ and $\gamma_k^b$.  In other words, if \eqref{eq:ODE-HB++} is to reproduce a heavy-ball method, the controller gains, the method of discretization and the discretization step-sizes must all be chosen jointly in some interdependent manner. Furthermore, even if it were somehow possible to choose the controller gains judiciously, generating a candidate algorithm by simply discretizing the resulting ODE using well-established numerical methods may not be prudent because ``the accuracy of the computed solution curve is not of prime importance''\cite{boggs71}; rather, it is more important to arrive at the ``asymptote of the solution ... with the fewest function evaluations''\cite{boggs71}.  In view of these observations, Boggs\cite{boggs71} proposed $A$-stable methods of integration to solve the differential equations that were previously generated by Davidenko and Gavurin\cite{gavurin}.  Despite his breakthrough, such methods are not widely used because they remain computationally expensive, a fact that has been known for quite sometime (see \cite{brown+biggs:sqp}). More recently, Grune and Karafyllis\cite{grune:optimization} developed a new idea based on framing a Runge-Kutta method as a hybrid dynamical system. In applying this approach to optimization, they concluded that ``if the emphasis lies on a numerically cheap computation ... then high order schemes may not necessarily be advantageous''\cite{grune:optimization}.

In pursuit of a new approach to generating algorithms, we choose to not produce the ODEs explicitly; instead, we revert back to the new foundations (cf.~\Cref{sec:new-principles}) that generated the ODEs in the first place.

%

\subsection{Development of a Three-Step Iterative Map}

In acknowledging that the needs of optimization are substantially different from those of traditional control theory as well as numerical methods for solving ODEs, we chart a new course for producing algorithms using the following ideas:
\begin{enumerate}
\item Rather than design the ODEs that generate the algorithm primitives, we directly use the minimum principles  within an algorithmic structure to find the instantaneous control $\bzeta(k)$ at iteration $k$.
\item Because an ODE that governs the algorithm primitive is never generated, we advance to the next iterate based on the geometric condition that every iterate remain on the zero-Hamiltonian singular manifold (cf.~\Cref{thm:H=0-infty}).
\end{enumerate}
The first idea leans on the concept of proximal aiming introduced by Clarke et al\cite{prox-aiming-1994} for an altogether different purpose of overcoming certain theoretical hurdles in nonsmooth control theory.  The second idea relies on using the readily available singular integral of motion (cf.~\eqref{eq:lamx=-gradL}) to generate $\blam_x(k+1)$ instead of discretizing and propagating its corresponding differential equation (cf.~\eqref{eq:lam_x-mod}).  Similarly, we generate $\bs(k+1)$ from \eqref{eq:idea-1} instead of discretizing \eqref{eq:idea-2}.  Consequently, only the simple linear equations in System~$(A)$ need be discretized.  Collecting all these ideas together, we arrive at the following procedure: Let $\zeta(k) := (\bu(k), \bmu(k), \omega(k))$ and $h_k > 0 $ be given.  Then a three-step iterative map for accelerated optimization is given by,
\begin{subequations}\label{eq:3-step-iMap}
\begin{align}
& A_1(k+1): \left\{
\begin{aligned}
\lambda_y(k+1) &=\lambda_y(k) + h_k\, \omega(k)\\
\blam_s(k+1) & = \blam_s(k) + h_k\, \bmu(k)\\
\bv(k+1) & = \bv(k) + h_k\, \bu(k)
\end{aligned}
\right. \\
& A_2(k+1): \left\{
\begin{aligned}
\quad \bx(k+1) &= \bx(k) + h_k \bv(k+1) \label{eq:BE4x}
\end{aligned}
\right. \\
& A_3(k+1): \left\{
\begin{aligned}
\blam_{x}(k+1)  &= -\partial_{\bx}L\big(\lambda_y(k+1), \blam_s(k+1), \bx(k+1)\big) \\
\bs(k+1) &=  \be\big(\bx(k+1)\big)
\end{aligned}
\right.
\end{align}
\end{subequations}
That is, $A_1(k+1)$ is used to generate the $(k+1)^{th}$ point for $\lambda_y, \blam_s$ and $\bv$ by a forward Euler method.  Next, the optimization variable $\bx$ is updated in $A_2(k+1)$ using a backward Euler formula.  Despite being a backward Euler formula, $A_2(k+1)$ is explicit because $\bv(k+1)$ is available from $A_1(k+1)$.  Finally, $\blam_x$ and $\bs$ are advanced to the $(k+1)^{th}$ point in $A_3(k+1)$ by not only sans discretization, but also that they are based on the most recent update of its arguments made available by $A_1(k+1)$ and $A_2(k+1)$.
\begin{remark}
Equation~\eqref{eq:3-step-iMap} is essentially a semi-discretization of System~$(A)$ (cf.~\eqref{eq:A4auxSystem}).
\end{remark}
\begin{remark}\label{rem:BE4x}
The backward Euler update for $\bx$ in \eqref{eq:BE4x} is essential not only for efficiency (i.e., in using the latest updates to generate new ones) but also to ensure consistency in terms of generating a Fritz John or KKT point.  This is because if a forward Euler method were to be used to update $\bx$ instead of \eqref{eq:BE4x}, then the sequence of iterates generated by \eqref{eq:3-step-iMap} will not advance to an improved point if $\bv_k$ were to vanish for some $k$ prior to achieving optimality.  Note also that \eqref{eq:BE4x} is implicitly contained in \eqref{eq:HB-gain-step}.
\end{remark}
%

From \eqref{eq:3-step-iMap} it follows that a feedback control law is never explicitly computed; hence, it is not necessary to produce an ODE that governs the flow of the algorithm primitive. Furthermore, because $\pounds_f V$ is linear in the control variable, Problems~$(P)$ and $(P^*)$ are ``simpler'' than the original problem $(N)$. In particular, note that Problem~$(P^*)$ is ``small'' scale; i.e., it has just one constraint equation, no matter the scale of the original problem $(N)$.

\subsection{Some New Step Length Procedures and Formulas}
As noted earlier, it is inadvisable to choose $h_k$ in \eqref{eq:3-step-iMap} based on the rules of numerical methods for ODEs. In view of this backdrop, we proposed in \cite{rossJCAM-1} a minimum principle for a maximal step length.  This principle is essentially an adaptation of the exact step length procedure used in standard optimization with the merit function replaced by the value of the CLF along the direction $\bzeta(k)$.  The key difference between the CLF and merit function approaches is that the former cannot be based on unconstrained optimization algorithms.  In advancing the maximal step-length principle for the iterative map given by \eqref{eq:3-step-iMap}, we pose the following problem for generating an exact step length $h_k$:
\begin{eqnarray}
& \bz: = (\blam_x, \lambda_y, \blam_s, \bv, \bs) \nonumber\\
&\qquad (P_h) \left\{
\begin{array} {llr}
\displaystyle\mathop\textsf{Minimize }_{\bz(k+1), \bx(k+1), h_k}  & V(\bz(k+1)) \\
\textsf{Subject to}     & \blam_{x}(k+1) + \partial_{\bx}L\big(\lambda_y(k+1), \blam_s(k+1), \bx(k+1)\big) = \bzero \\
                        & \lambda_y(k+1) - \lambda_y(k) - h_k\, \omega(k) = 0\\
                        & \blam_s(k+1) - \blam_s(k) - h_k\, \bmu(k) = \bzero\\
                        & \bv(k+1) - \bv(k) - h_k\, \bu(k) = \bzero \\
                        & \bs(k+1) - \be\big(\bx(k+1)\big) = \bzero \\
                        & \bx(k+1) - \bx(k) - h_k \bv(k+1) = \bzero \\
                        &  h_k  \ge 0
\end{array} \right.& \label{prob:max-h}
\end{eqnarray}
Assuming $h_k > 0$, the dual feasibility conditions for Problem~$(P_h)$ are given by,
\begin{subequations}\label{eq:dualFeas}
\begin{align}
\partial_{\bz} V(\bz(k+1)) + \left[
                               \begin{array}{c}
                                 \bpsi_{\lambda_x} \\
                                 \psi_{\lambda_y} + \partial_{\bx}E(\bx(k+1))\cdot \bpsi_{\lambda_x} \\
                                 \bpsi_{\lambda_s} + [\partial_{\bx}\be(\bx(k+1))] \bpsi_{\lambda_x} \\
                                 \bpsi_v - h_k \bpsi_x \\
                                 \bpsi_s \\
                               \end{array}
                             \right]  & = \bzero\\
[\partial^2_{\bx}L(\lambda_y(k+1), \blam_s(k+1), \bx(k+1))]\bpsi_{\lambda_x} - \partial_{\bx}\be(\bx(k+1)) \cdot \bpsi_{s} + \bpsi_x  &= \bzero \\
\psi_{\lambda_y} \omega(k)  + \bpsi_{\lambda_s} \cdot \bmu(k) + \bpsi_v \cdot \bu(k) - \bpsi_x \cdot \bv(k+1) & = 0 \label{eq:dualFeas-hk}
\end{align}
\end{subequations}
where, $\bpsi_{\lambda_x}, \psi_{\lambda_y}, \bpsi_{\lambda_s}, \bpsi_v, \bpsi_s$ and $\bpsi_x$ are Lagrange multipliers associated with the constraint equations given in \eqref{prob:max-h}.
\begin{remark}
Because Problem~$(P_h)$ incorporates \eqref{eq:3-step-iMap}, it also generates the iterates to solve Problem~$(N)$; i.e., if Problem~$(P_h)$ can be solved ``exactly,'' then its solution, together with that of any one of the minimum principles represents the complete algorithm.
\end{remark}

It is apparent by a cursory inspection of the primal and dual feasibility conditions of Problem~$(P_h)$ that producing an explicit equation for $h_k$ in terms of the known information at $k$ is not readily possible; in fact, this challenge is not altogether different than the problem of generating an exact step length formula using standard merit functions.  In view of this, it is apparent that $h_k$ may be generated more efficiently by using the traditional approach of inexact line search methods (i.e., Armijo-Goldstein-Wolfe methods) but adapted to the values of the CLF along the direction $\bzeta(k)$.  Nonetheless, as is well-known, the efficiency of such methods are more strongly dependent on the initial value of $h_k$ rather than the specifics of backtracking.  Consequently, motivated by the need to produce a ``good'' initial value for $h_k$, we advance three formulas.

The first approach is based on approximating the constraints in Problem~$(P_h)$ and solving the resulting problem.  The constraint approximations are based on the first order terms in $h_k$; this generates the following approximations:
\begin{subequations}\label{eq:approx4hk}
\begin{align}
&\partial_{\bx} L(\lambda_y(k+1), \blam_s(k+1), \bx(k+1)) \approx \partial_{\bx} L(\lambda_y(k), \blam_s(k), \bx(k)) \nonumber \\
&\qquad + h_k\, \partial^2_{\bx} L(\lambda_y(k), \blam_s(k), \bx(k))\, \bv(k+1) + h_k \, \partial_{\bx} L(\omega(k), \bmu(k), \bx(k)) \\
& \be(\bx(k+1)) \approx \be(\bx(k)) + h_k\, \partial_{\bx}\be(\bx(k))  \bv(k+1) \\
& \partial^2_{\bx} L(\lambda_y(k+1), \blam_s(k+1), \bx(k+1)) \approx \partial^2_{\bx} L(\lambda_y(k), \blam_s(k), \bx(k))\\
&\partial_{\bx}\be(\bx(k+1)) \approx \partial_{\bx}\be(\bx(k)) \\
& \partial_{\bx}E(\bx(k+1)) \approx \partial_{\bx}E(\bx(k))
\end{align}
\end{subequations}
%
\begin{proposition}
Suppose a CLF is given by the quadratic function $V(\bz) = (\bz^T\bQ\bz)/2$ where, $\bQ$ is a positive definite matrix. Assume \eqref{eq:approx4hk} holds.  Let $\bz_a(k+1)$ denote the approximate value of $\bz$ based on the approximations given by \eqref{eq:approx4hk}.   Then, a solution to $h_k = h^{BL}_k$ satisfies the system of bilinear equations given by,
\begin{equation}\label{eq:bilinear4hk}
\bM_1(k) \left[
      \begin{array}{c}
        \bchi \\
        h_k^{BL} \\
      \end{array}
    \right] + h_k^{BL} \bM_2(k)\, \bchi = \bb_k
\end{equation}
where, $\bM_1(k), \bM_2(k)$ and $ \bb_k$ are matrices (of appropriate dimensions) that depend on the known values of the iterates of \eqref{eq:3-step-iMap} at the point $k$, and $\bchi$ is a variable that comprises $\bz_a(k+1)$, $\bpsi_{\lambda_x}, \psi_{\lambda_y}, \bpsi_{\lambda_s}, \bpsi_v, \bpsi_s$ and $\bpsi_x$.
\end{proposition}
\begin{proof}
The result follows from two simple steps:  Replace $\bv(k+1)$ in \eqref{eq:dualFeas-hk} by $\bv(k) + h_k \bu(k)$. Substitute \eqref{eq:approx4hk} in \eqref{eq:dualFeas} and \eqref{eq:3-step-iMap}.
\end{proof}
\begin{remark}
The approximations given by \eqref{eq:approx4hk} are only used to generate $h_k^{BL}$ via \eqref{eq:bilinear4hk}.  In other words, $\bz_a(k+1)$ generated from \eqref{eq:bilinear4hk} is discarded once $h_k^{BL}$ is computed.
\end{remark}

A second formula for an initial value of $h_k$ is given by the following proposition:

\begin{proposition}[\cite{rossJCAM-1}]
Let $V(\bz) = (\bz^T\bQ\bz)/2$ where, $\bQ$ is a positive definite matrix.  Suppose all of the constraint equations in \eqref{prob:max-h} are replaced by a forward Euler discretization. Then, the resulting maximal step length is given explicitly by,
\begin{equation}\label{eq:h-VquadEuler-explicit}
h^{FE}_k = -\frac{\bz_k^T\bQ\bff_k}{\bff_k^T\bQ\bff_k} = \frac{-\pounds_f V(\bz_k)}{2 V(\bff_k)}
\end{equation}
where $\bff_k = \bff(\bz_k, \bzeta_k)$ and $\bff$ is given by \eqref{eq:f=sumOf2}.
\end{proposition}

Because the minimum principles $(P)$ and $(P^*)$ generate $\pounds_f V(\bz_k) < 0$, it follows that \eqref{eq:h-VquadEuler-explicit} guarantees $h^{FE}_k > 0$.   The main problem in using $h^{FE}_k$ in \eqref{eq:3-step-iMap} is its inconsistency as discussed in \Cref{rem:BE4x}; however, it holds the potential of providing a lower bound for an acceptable step size in a Goldstein-type condition.

Finally, a third formula for an initial value of $h_k$ is obtained by utilizing the fact that $\pounds_f V(\bz_k)$ is the continuous-time derivative of $V$ at the point $\bz_k$. As a result, the tangent line emanating from the point $\bz_k$ may be parameterized as,
\begin{equation}\label{eq:V-tangent-line}
V^{tan}(s) = s \pounds_f V(\bz_k)  + V(\bz_k)
\end{equation}
Setting $V^{tan}(s) = 0$ in \eqref{eq:V-tangent-line} to solve for $s$ as a proposed value for an initial step size generates the very simple formula,
\begin{equation}\label{eq:h=h-tan}
h^{tan}_k = \frac{V(\bz_k)}{-\pounds_f V(\bz_k)}
\end{equation}
%

\subsection{Description of the Main Algorithm}
The main algorithm  comprises two key steps:
\begin{enumerate}
\item At step $k$, solve Problem~$(P)$ or $(P^*)$ to generate $\bzeta(k)$.  For a quadratic CLF, this only requires a solution to a linear system; see \eqref{eq:MinP4search} and \eqref{eq:MinP*4search}.
\item Using the computed value of $\bzeta(k)$ from the prior step, advance to step $(k+1)$ using \eqref{eq:3-step-iMap} such that $V_{k+1}$ is sufficiently less than $V_k$, where $V_{k+1}$ is the value of the CLF at the accepted point $(k+1)$.
\end{enumerate}
The major steps of the proposed algorithm are encapsulated in \cref{alg:main}.
\begin{algorithm}
\caption{Main}
\label{alg:main}
\begin{algorithmic}
\STATE{ Choose a CLF $V$ and the parameters associated with Problem~$(P)$ or $(P^*)$ (cf.~\eqref{eq:MinP4search} and \eqref{eq:MinP*4search})}
\STATE{ Initialize the algorithm according to: $\blam_x^0 = -\partial_{\bx}L(1, \blam_s^0, \bx^0)$, $\bs^0 = \be(\bx^0)$. 
Set $k = 0$. Compute $V_k = V(\bz(k))$}
\WHILE{stopping conditions are not met}
\STATE{Generate $\bzeta(k)$ by solving Problem~$(P^*)$ (or $(P)$)}
\STATE{Compute $h_k^0$ using any one of \eqref{eq:bilinear4hk}, \eqref{eq:h-VquadEuler-explicit} or \eqref{eq:h=h-tan} }
\STATE{Advance to $\big(\bz(k+1), \bx(k+1)\big)$ using \eqref{eq:3-step-iMap} and $h_k^0$}
\STATE{Compute $V_{k+1} = V(\bz(k+1))$}
\WHILE{$V_{k+1}$ has not decreased sufficiently}
\STATE{Backtrack $\big(\bz(k+1), \bx(k+1)\big)$ along $\bzeta(k)$; recompute $V_{k+1}$}
\ENDWHILE
\STATE{Update $k \leftarrow  k + 1$}
\ENDWHILE
\end{algorithmic}
\end{algorithm}

\subsection{A Numerical Illustration} \label{sec:example}
\begin{figure}[h]
  \centering
    \begin{subfigure}{0.45\textwidth}
        \centering
        \includegraphics[width=0.9\linewidth]{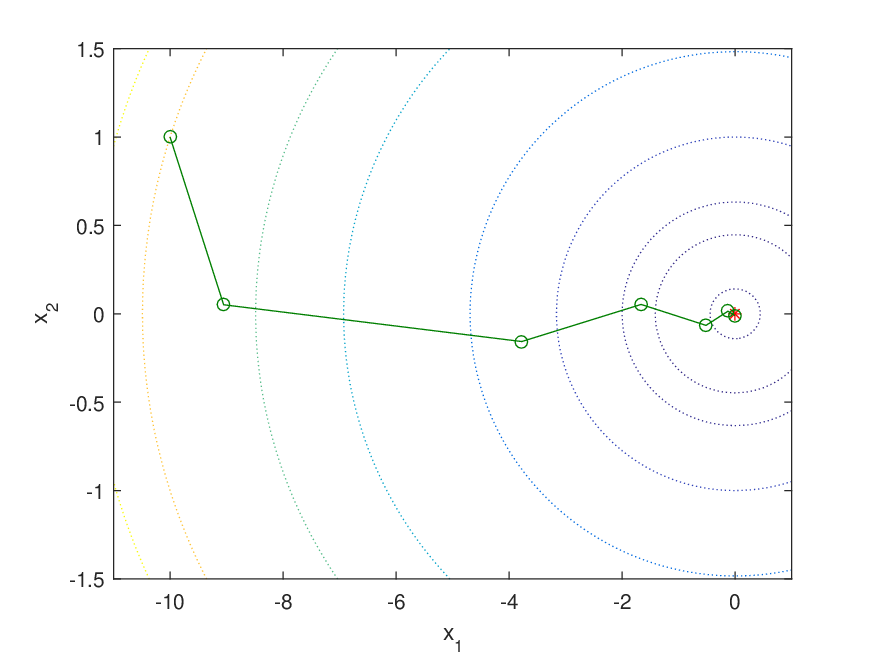}
        \caption{Six iterations of a gradient algorithm based on \Cref{alg:main}}
         \label{fig:a}
        \end{subfigure}
\hfill
    \begin{subfigure}{0.45\textwidth}
        \centering
        \includegraphics[width=0.9\linewidth]{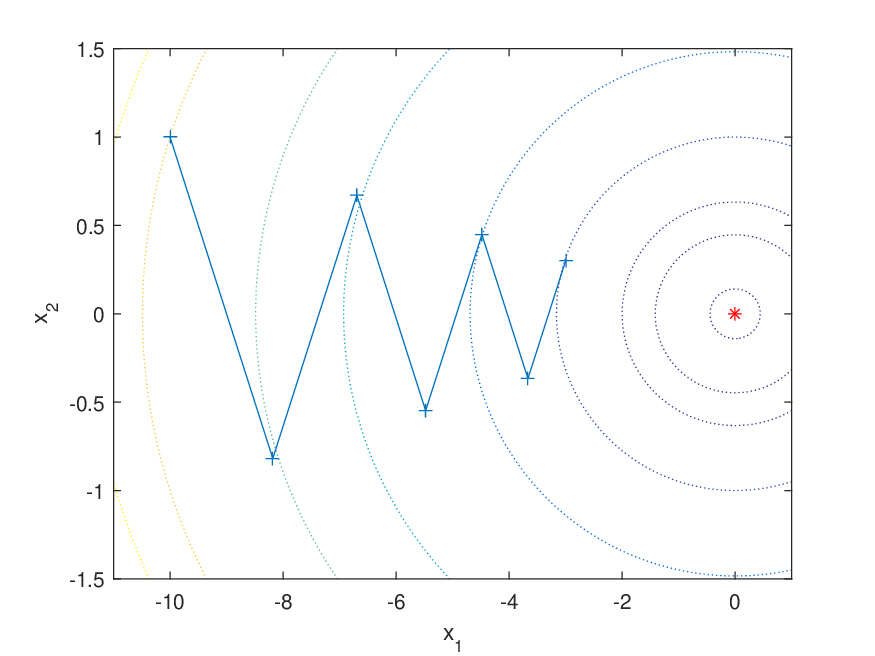}
        \caption{Six iterations of a standard gradient algorithm}
        \label{fig:b}
        \end{subfigure}
\newline
 \begin{subfigure}{0.45\textwidth}
        \centering
        \includegraphics[width=0.9\linewidth]{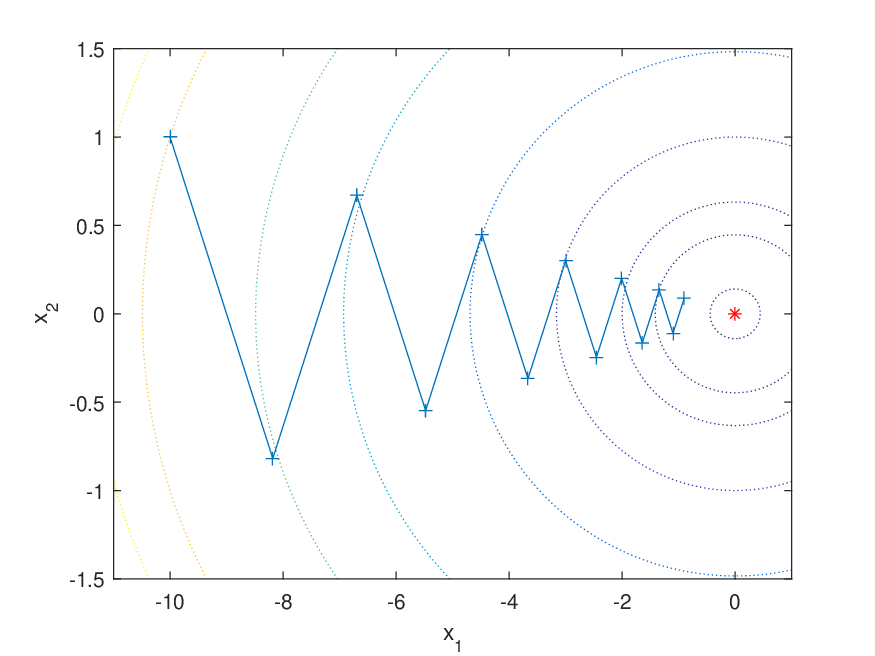}
        \caption{Twelve iterations of a standard gradient algorithm}
         \label{fig:c}
        \end{subfigure}
  \hfill
    \begin{subfigure}{0.45\textwidth}
        \centering
        \includegraphics[width=0.9\linewidth]{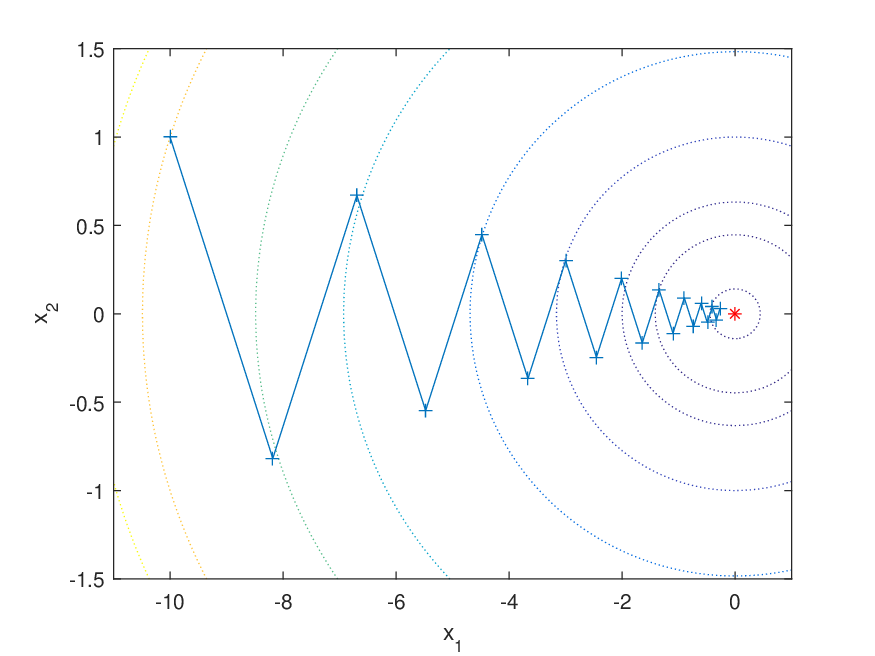}
        \caption{Eighteen iterations of a standard gradient algorithm}
        \label{fig:d}
        \end{subfigure}
  \caption{An illustration of an application of \Cref{alg:main} for accelerating a gradient method.}
  \label{fig:quadPlots}
\end{figure}

We present a simple numerical example to demonstrate the acceleration generated by an application of \Cref{alg:main}.  Shown in \Cref{fig:a} are six iterates of \Cref{alg:main} applied to minimize the function $(x_1, x_2) \mapsto (x_1^2 + 10x_2^2)/2$.  The iterates were obtained by setting $\bW$ to be identity matrix (cf.~\Cref{sec:derive-HB+}); hence the resulting algorithm is a new gradient method. To demonstrate the fact that this new gradient method does indeed achieve acceleration, the iterates of a standard gradient method for the same number of iterations (i.e., six) are shown in \Cref{fig:b}.
To provide an additional perspective in terms of an acceleration factor generated by an application of \Cref{alg:main}, iterates of the standard gradient method for double and triple the number of iterations are shown in \Cref{fig:c} and \Cref{fig:d} respectively.


\section{Conclusions}\label{sec:conclusions}
The transversality mapping principle and its consequences facilitate new ideas for designing and analyzing optimization algorithms. A general framework for accelerated (and unaccelerated) optimization methods is possible under the rubric singular optimal control theory.  On hindsight, the central role of singular optimal control theory is not surprising because a nonsingular control would have implied a universal optimal algorithm.  By the same token, the infinite-order of the singular arc is also not surprising because a finite order would also imply a universal optimal algorithm.  The interesting aspect of many well-known algorithms -- accelerated or otherwise -- is that their primitives are all describable in terms of flows over a zero-Hamiltonian singular manifold. This insight is used to launch a three-step iterative map that generates iterates which remain on the singular manifold. It turns out that the key steps to computational efficiency is not necessarily based on discretizing the resulting ordinary differential equations, rather, it is based on combining the more traditional aspects of optimization with the generation of Euler polygonal arcs by proximal aiming.  There is no doubt that a vast number of open questions remain; however, it is evident that new viable optimization algorithms can indeed be generated using the results emanating from the transversality mapping principle.

%

\bibliographystyle{siamplain}

\begin{thebibliography}{10}

\bibitem{rossJCAM-1}
I. M. Ross, An optimal control theory for nonlinear optimization, J. Comp. and Appl. Math., 354 (2019) 39--51.

\bibitem{vinter}
R. B. Vinter,  Optimal Control, Birkh\"{a}user, Boston, MA,
2000.

\bibitem{ross-book}
I. M. Ross, A Primer on Pontryagin's Principle in Optimal Control, second ed., Collegiate Publishers, San Francisco, CA, 2015.

\bibitem{krener-hmp}
A. J. Krener, The high order Maximal Principle and its applications to singular extremals, SIAM J. of control and optimization, 15/2 (1977), 256--293.

\bibitem{clarkeLyap}
F. Clarke, Lyapunov functions and feedback in nonlinear control. In: M.S. de Queiroz, M. Malisoff, P. Wolenski (eds) Optimal control, stabilization and nonsmooth analysis. Lecture Notes in Control and Information Science, vol 301. Springer, Berlin, Heidelberg (2004), 267--282.

\bibitem{clarkeEncy}
F. Clarke, Nonsmooth analysis in systems and control theory. In: Meyers, R. A. (ed) Encyclopedia of Complexity and Systems Science. Springer, New York, N.Y. (2009), 6271--6285.

\bibitem{sontag-book}
E. D. Sontag, Mathematical Control Theory: Deterministic Finite Dimensional Systems, second ed., Springer, New York, NY, 1998.

\bibitem{motta-CLFs}
M. Motta, F. Rampazzo, Asymptotic controllability and Lyapunov-like functions determined by Lie brackets, SIAM J. Control and Optimization, 56/2, 2018, pp.~1508--1534.


\bibitem{freeman-acc}
R. A. Freeman, P. V. Kokotovi\'{c}, Optimal nonlinear controllers for feedback linearizable systems, Proc. ACC, Seattle, WA, June 1995.

\bibitem{bhat-2000}
S. P. Bhat, D. S. Bernstein, Finite-time stability of continuous autonomous systems, SIAM J. Control Optim., 38/3, 2000, pp.~751--766.

\bibitem{CLFtoHJB2020}
P. Osinenko, P. Schmidt, S. Streif, Nonsmooth stabilization and its computational aspects, IFAC-PapersOnLine, 53/2, 2020, pp.~6370--6377,

\bibitem{NLP2ODE-1980}
H. Yamashita, A differential equation approach to nonlinear programming, Mathematical Programming, 18, 1980, pp.~155--168.

\bibitem{NLP2ODE-1981}
D. M. Murray, S. J. Yakowitz, The application of optimal control methodology
to nonlinear programming problems. Math. Programming, 21/3, 1981, pp.~331--347.

\bibitem{NLP2ODE-1989}
A. A. Brown, M. C. Bartholomew-Biggs, ODE versus SQP methods for constrained optimization, J. optimization theory and applications, 62/3, 1989, pp.~371--386.

\bibitem{NLP2ODE-1994}
Yu.G. Evtushenko, V.G. Zhadan, Stable barrier-projection and barrier-Newton methods in nonlinear programming, Optim. Methods Software 3, 1994, pp.~237--256.

\bibitem{NLP2ODE-2006}
A. Bhaya, E. Kaszkurewicz, Control Perspectives on Numerical Algorithms and Matrix Problems, Advances in Design and Control, SIAM, Philadelphia, PA, 2006.

\bibitem{NLP2ODE-2007}
L. Zhou, Y. Wu, L. Zhang, and G. Zhang, Convergence analysis of a differential equation
approach for solving nonlinear programming problems, Appl. Math. Comput., 184, 2007,
pp.~789--797.

\bibitem{NLP2ODE-2017}
I. Karafyllis, M. Krstic, Global Dynamical Solvers for Nonlinear Programming Problems, SIAM J. Control and Optimization, 55/2, 2017, pp.~1302--1331.


\bibitem{lessard-2016}
L. Lessard, B. Recht, A. Packard, Analysis and design of optimization algorithms via integral quadratic constraints, SIAM Journal on Optimization, 2016, 26(1), 57--95.

\bibitem{su-2016}
W. Su, S. Boyd, E. J. Candes, A differential equation for modeling Nesterov's accelerated gradient method: theory and insights, J. machine learning research, 17 (2016) 1--43.

\bibitem{wibisono-2016}
A. Wibisono, A. C. Wilson,  M. I. Jordan, A variational perspective on
accelerated methods in optimization, Proceedings of the National Academy of Sciences, 2016,
133:E7351--E7358.


\bibitem{Goh-1997}
B. S. Goh, Algorithms for unconstrained optimization via control theory, J. Optim. Theory Appl., 92/3, 1997, pp.~581--604.

\bibitem{Goh-2021}
M. S. Lee, H. G. Harno, B. S. Goh, K. H. Lim, On the bang-bang control approach via a component-wise line search strategy for unconstrained optimization, Numerical Algebra, Control and Optimization, 11/1, 2021, pp.~45--61.

\bibitem{polyak64}
B. T. Polyak, Some methods of speeding up the convergence of iteration methods, USSR Computational Math. and Math. Phys., 4/5 (1964) 1--17 (Translated by H. F. Cleaves).

\bibitem{polyak2017}
B. Polyak, P. Shcherbakov, Lyapunov functions: an optimization theory perspective, IFAC PapersOnLine, 50-1 (2017) 7456--7461.

\bibitem{nesterov83}
Yu. E. Nesterov, A method of solving a convex programming problem with convergence rate $\mathcal{O}(1/k^2)$, Soviet Math. Dokl., 27/2 (1983) 371--376 (Translated by A. Rosa).






\bibitem{boggs71}
P. T. Boggs, The solution of nonlinear system of equations by $A$-stable integration techniques, SIAM J. Numer. Anal. 8/4 (1971) 767--785.

\bibitem{gavurin}
M. K. Gavurin, Nonlinear functional equations and continuous analogues of iteration methods, Izv. Vyssh. Uchebn. Zaved. Mat., 5 (1958) 18--31.

\bibitem{brown+biggs:sqp}
A. A. Brown, M. C. Bartholomew-Biggs, ODE versus SQP methods for constrained optimization, J. optimization theory and applications, 62/3 (1989) 371--386.

\bibitem{grune:optimization}
L. Gr\"{u}ne,  I.  Karafyllis,  Lyapunov Function Based Step Size Control for Numerical
ODE  Solvers  with  Application  to  Optimization  Algorithms. In: K. H\"{u}per, J. Trumpf (eds.)
Mathematical System Theory, pp.~183--210 (2013) Festschrift in Honor of Uwe Helmke on the Occasion of his 60th Birthday.


\bibitem{prox-aiming-1994}
F. H. Clarke, Yu. S. Ledyaev, and A. I. Subbotin. Universal feedback control via proximal aiming in problems of control under disturbances and differential games. Univ. de Montr\'{e}al, Report CRM 2386, 1994.

\end{thebibliography}

\end{document}